\documentclass[12pt]{amsart}
\usepackage{amsmath}
\usepackage{amsfonts}
\usepackage{amsfonts}
\usepackage{amssymb}
\usepackage{amsthm}
\usepackage{amsmath}
\usepackage{amstext}
\usepackage[T2A]{fontenc}
\usepackage[cp1251]{inputenc}

\usepackage{amsmath}
\usepackage{amsmath, pb-diagram}
\usepackage{pb-diagram}
\usepackage[all]{xy}
\usepackage{amssymb}
\usepackage{amscd}
\newtheorem{thm}{Theorem}[section]
 \newtheorem{cor}[thm]{Corollary}
 \newtheorem{lem}[thm]{Lemma}
 \newtheorem{prop}[thm]{Proposition}

\theoremstyle{definition}
\newtheorem {ex}[thm] {Example}
\newtheorem {rem} [thm] {~~~Remark}

\theoremstyle{plain}

\newcommand{\End}{\operatorname{End}}

\newcommand{\Soc}{\operatorname{Soc}}

\newcommand{\Aut}{\operatorname{Aut}}

 \numberwithin{equation}{section}

\DeclareMathOperator{\Z}{\mathbb{Z}}

\begin{document}

\author{Adel  Abyzov,  Truong Cong Quynh, Askar Tuganbaev}

\address{ Department of Algebra and Mathematical Logic, Kazan (Volga Region) Federal University, 18 Kremlyovskaya str., Kazan, 420008 Russia}
\email{ aabyzov@ksu.ru}

\address{Department of Mathematics, The University of Danang - University of Science and Education,  459 Ton Duc Thang, Danang city, Vietnam}
\email{tcquynh@ued.udn.vn}

\address{National Research University ``Moscow Power Engineering Institute''}
\email{tuganbaev@gmail.com}

\title[ Rings whose ideals are close to automorphism-invariant
]
 {Rings whose ideals are close to automorphism-invariant
}
\maketitle

\maketitle
\begin{abstract}  We consider rings whose one-sided ideals are close to automorphism-invariant modules. We study rings in which every (finitely generated) right ideal is automorphism invariant and rings in which every right ideal is a finite direct sum of automorphism invariant ideals. Connections between these classes of rings, $q$-ring and $\Sigma$-$q$-rings are also considered.
\end{abstract}
\vskip 2cm

\section{Introduction}

All rings are assumed associative and unitary, and all modules are unitary.
A module is quasi-injective provided that it is invariant under the endomorphisms of its injective
envelope. A module is automorphism-invariant provided that it is invariant under the automorphisms
of its injective envelope. Automorphism-invariant modules were first considered by Dickson and Fuller
in \cite{DF} and were systematically studied  recently.
Characterization of rings by the properties of their one-sided ideals is one of the important areas
in ring and module theories. A ring $R$ is a right $a$-ring ($q$-ring, respectively) provided that its every
right ideal is an automorphism-invariant right $R$-module (quasi-injective, respectively). The $q$-rings were
studied in many articles (for instance, see \cite{BFKJ}, \cite{H73}, \cite{H73}, \cite{I72}, \cite{I75}, \cite{JMS}, \cite{B79}). Full description of these rings was obtained in \cite{BFKJ}.
Various properties of $a$-rings were studied in \cite{TQS}. In particular, as shown in \cite{TQS}, a right $a$-ring is regular
if and only if it is semiprime, and the structure of the right Artinian nonsingular right $a$-rings was found. In \cite{I96}, Ivanov generalized the $q$-rings, namely $fq$-rings,
that the rings all of whose finitely generated (right) ideals are quasi-injective. He
studied $fq$-rings in conjunction with dense primitive idempotents and idempotentnonsingular rings and obtained some interesting results. He showed that the main
structure theorem on injective von Neumann regular rings can be extended to
$fq$-rings which are idempotent-nonsingular and determined the structure of indecomposable idempotent-nonsingular fq-rings with dense primitive idempotents and
represented them as rings of matrices.  Ring $R$ is called a  right $fa$-ring  if every finitely generated right ideal of $R$ is  automorphism-invariant. Various properties of $fa$-rings were studied in the article \cite{QAT}.

A ring $R$ is a right $\Sigma$-$a$-ring ($\Sigma$-$q$-ring, respectively) provided that its every right ideal is a finite
direct sum of automorphism-invariant right $R$-modules (quasi-injective, respectively). The $\Sigma$-$q$-rings were
introduced and studied in \cite{JSS}. Some properties of $q$-rings and $\Sigma$-$q$-rings were reflected in \cite{JST12} as well.
The $\Sigma$-$a$-rings were first introduced in \cite{SS2}, and the question of description of $\Sigma$-$a$-rings was posed there
(Question 4, p. 310). Various properties of $\Sigma$-$a$--rings were studied in \cite{APT}. In particular, it was established a representation in the form of block
upper triangular rings of formal matrices for the indecomposable right Artinian right hereditary right
$\Sigma$-$a$-rings.

In the present article, provides an overview of recent results related to $a$-rings, $fa$-rings and $\Sigma$-$a$-rings.

The Jacobson radical (the right singular ideal, respectively) of a ring R is denoted by $J(R)$ (by $Z_r(R),$
respectively). The fact that $N$ is a submodule of a module $M$ is denoted by $N\leq  M.$ The injective
envelope of a module $M$ is denoted by $E(M).$ If $M$ is a module of finite length then its length is denoted
by $l(M).$

\section{Rings with every right ideal automorphism-invariant.}

In 1969, Jain and  Singh in \cite{7} introduced rings for which every right ideal is quasi-injective, it is called {\it right q-rings}. Rings all of whose right ideals are automorphism-invariant are called {\it right $a$-rings} see  \cite{SS2} or \cite{TQS}.  Recall that right $q$-rings are precisely those right self-injective rings for which every essential right ideal is a two-sided ideal. In particular,   commutative self-injective rings are   $q$-rings and hence  $a$-rings. Next we would like to present some examples of   $a$-rings that are not $q$-rings.

\begin{ex}
	Consider the ring $R$ consisting of all eventually constant sequences of elements from $\mathbb F_2$. Clearly, $R$ is a commutative automorphism-invariant ring as the only automorphism of its injective envelope is the identity automorphism. Hence $R$ is an $a$-ring by the above lemma. But $R$ is not a $q$-ring because $R$ is not self-injective.
\end{ex}

Next,  we will consider some properties for right $a$-rings. These equivalent characterizations will be more convenient to use.

\begin{prop}{\cite[Proposition 3.1]{TQS}}\label{essen} The following conditions are equivalent for a ring $R$:
	\begin{enumerate}
		\item $R$ is a right $a$-ring.
		\item Every essential right ideal of $R$ is automorphism-invariant.
		\item $R$ is right automorphism-invariant and every essential right ideal of $R$ is a left $T$-module, where  $T$ is a subring of $R$ generated by its unit elements.
	\end{enumerate}
\end{prop}

We can show that for a right $a$-ring and $A$ be a right ideal of $R$, if there exists a right ideal $B$ of $R$ with $A\cap B=0$ and $A\cong B$, then   $A$ is nonsingular semisimple and injective.

\noindent Recall that two modules $M$ and $N$ are said to be {\it orthogonal} if no submodule of $M$ is isomorphic to a submodule of $N$. A module $M$ is said to be a {\it square module} if there exists a right module $N$ such that  $M\cong N^2$ and a submodule $N$ of a module $M$ is called {\it square-root} in $M$ if $N^2$ can be embedded in $M$. A module $M$ is called {\it square-free} if $M$ contains no non-zero  square roots and $M$ is called {\it square-full} if every submodule of $M$ contains a non-zero square root in $M$.

\begin{thm}{\cite[Theorem  3.4]{TQS}}\label{thm:g} A right $a$-ring is a direct sum of a square-full semisimple artinian
	ring and a right square-free ring.
\end{thm}

As a consequence of the above, we have

\begin{cor}
	An indecomposable ring $R$ containing a square is a right $a$-ring if and only if $R$ is simple artinian.
\end{cor}

We denote the ring of $n\times n$ matrices over a ring $R$ by $\mathbb M_n(R)$. In the next theorem we study when matrix rings are right $a$-rings.

\begin{thm}{\cite[Theorem  3.6]{TQS}}\label{pro:ma} Let $n>1$ be an integer. The following conditions are equivalent for a ring $R$:
	\begin{enumerate}
		\item $\mathbb{M}_n(R)$ is a right $q$-ring.
		\item $\mathbb{M}_n(R)$ is a right $a$-ring.
		\item $R$ is semisimple artinian.
	\end{enumerate}
\end{thm}
 
The following example shows that there exists automorphism-invariant rings which are not right $a$-rings.

\begin{ex}Let $R=\Z_{p^n}$, where $p$ is a prime and $n>1$. It is well known that  $R$ is self-injective. Then, $\mathbb{M}_m(R)$ is right self-injective for all $m> 1$. Thus, for instance, $\mathbb{M}_m(\Z_{p^2})$ is a right automorphism-invariant ring. But $\mathbb{M}_m(\Z_{p^2})$ is not a right $a$-ring for any $m>1$ in the view of above theorem as $\Z_{p^2}$ is not semisimple artinian. This example also shows that being a right $a$-ring is not a Morita invariant property.
\end{ex}

We continue to  consider some special classes of rings, for example, simple, semiprime, prime   and characterize as to when these rings are right $a$-rings. We begin this part with a simple observation.

\noindent Recall that a ring $R$ is called {\it von Neumann regular} if for every $a\in R$, there
exists some $b\in R$ such that $a = aba$. A ring $R$ is said to be {\it prime} if the product of any two nonzero ideals of $R$ is nonzero and a ring $R$ is called {\it semiprime} if it has
no nonzero nilpotent ideals.

\begin{thm}{\cite[Theorem  4.2]{TQS}} \label{reg} A right $a$-ring is von Neumann regular if and only if it is semiprime.
\end{thm}
 
A ring $R$ is called {\it unit-regular} if, for every element $x\in R$, there exists a unit $u \in  R$ such that $x = xux$. We can now have the following result.

\begin{cor} \label{vnr}
	Every von Neumann regular right $a$-ring is unit-regular.
\end{cor}

\begin{cor} The ring of linear transformations $R:=\End(V_D)$ of a vector space $V$  over a division ring $D$  is a right $a$-ring if and only if the vector space is finite-dimensional.
\end{cor}

\noindent A ring $R$ is said to be {\it strongly regular} if for every $a\in R$,
there exists some $b\in R$ such that $a = a^2b$.

\begin{prop}{\cite[Proposition 4.6]{TQS}} Let $R$ be a semi-prime right $a$-ring with zero socle. Then $R$ is strongly
	regular.
\end{prop}

\begin{thm} {\cite[Theorem  4.7]{TQS}}\label{prime}
	Let $R$ be a prime ring. Then $R$ is a right $a$-ring if and only if $R$ is a simple artinian ring.
\end{thm}

In particular, from the above theorem it follows that every simple right $a$-ring is artinian.

Next,   we would like to describe the structure of right $a$-rings.  
\begin{thm} $($Beidar, Fong, Ke, Jain, \cite{BFKJ}$)$
	A right $q$-ring $R$ is isomorphic to a finite direct product of right $q$-rings of the following types:
	\begin{enumerate}
		\item Semisimple artinian ring.
		\item $H(n; D; id_D)$ where $id_D$ is the identity automorphism on division ring $D$.
		\item $G(n; \Delta; P)$ where $\Delta$ is a right $q$-ring whose all
		idempotents are central.
		\item A right $q$-ring whose all idempotents are central.
		
	\end{enumerate}
	
	Here\\

	$H(n; D; \alpha)= \left[
	\begin{array}{ccccccc}
		D & V & 0 & &  &  & 0 \\
		0 & D & V & 0 & &  & 0 \\
		&  & D & V & 0 &  &  \\
		&  &  &  &  &  &  \\
		&  &  &  & D & V & 0 \\
		&  &  &  &  & D & V \\
		V(\alpha ) & 0 & &  &  &  & D%
	\end{array}
	\right] $, where $V$ is one-dimensional both as a left $D$-space and a right
	$D$-space, $V(\alpha )$ is also a one-dimensional left $D$-space as well as
	a right $D$-space with right scalar multiplication twisted by an
	automorphism $\alpha $ of $D$, i.e., $vd = v\cdot\alpha(d)$ for all $v\in V$, $d\in D$,
	\newline and
	$$
	G_n(n; \Delta; P):=\begin{pmatrix}
		D&V&&&&  \\
		&D&V&&&\\
		&&D&V&&\\
		&&.&.&.&\\
		&&&&.&.&\\
		&&&&.&D&V\\
		&&&&&&\Delta\\
	\end{pmatrix},$$
	where $V$ is as above and $\Delta$ is a right $q$-ring with maximal essential right
	ideal $P$ and hence $D=\Delta/P$ is a division ring.
\end{thm}

Now, using the above defined notations, we give the following description of right $a$-rings.

\begin{thm}{\cite[Theorem  5.2]{TQS}} \label{ornn} Let $n\ge 1$ be an integer, $D_1, D_2,\dots,D_n$ be division rings and $\Delta$ be a right $a$-ring with all idempotents central and an essential
	ideal, say $P$, such that $\Delta/P$ is a division ring and the right $\Delta$-module $\Delta/P$ is not embeddable into $\Delta_{\Delta}$. Next, let $V_i$ be a $D_i$-$D_{i+1}$-bimodule such that
	$$dim ({}_{D_i}\{V_i\})=dim(\{V_i\}_{D_{i+1}})=1$$ for all $i=1,2,\dots,n-1$,  and let $V_n$ be a $D_n$-$\Delta$-bimodule such that $V_nP=0$ and $$dim ({}_{D_n}\{V_n\})=dim(\{V_n\}_{\Delta/P})=1.$$
	Then $R:= G_n(D_1,\dots,D_n, \Delta, V_1,\dots, V_n)$ is a right $a$-ring.
	\end {thm}

	We finish this section by giving by another structure theorem for indecomposable right artinian right non-singular right $a$-ring describing them as a triangular matrix ring of certain block matrices.

	A ring  $R$ is called a {\it right  weak CS ring} if every semisimple right ideal  of $R$ is essential in a direct summand of $R_R$.

	\begin{thm}{\cite[Theorem  5.6]{TQS}} \label{ornn1}
		Any indecomposable right artinian right nonsingular right weakly CS right $a$-ring $R$ is isomorphic to
		$$ \begin{pmatrix} \mathbb{M}_{n_1}(e_1Re_1)&\mathbb{M}_{n_1\times n_2}(e_1Re_2)&\mathbb{M}_{n_1\times n_3}(e_1Re_3)&\cdots&\mathbb{M}_{n_1\times n_k}(e_1Re_k)\\
			0&\mathbb{M}_{n_2}(e_2Re_2)&\mathbb{M}_{n_2\times n_3}(e_1Re_2)&\cdots &\mathbb{M}_{n_2\times n_k}(e_2Re_k)\\
			0&0&.& \cdots & .\\
			\vdots &\vdots& \vdots & \vdots & \vdots \\
			0&0&0&\cdots&\mathbb{M}_{n_k}(e_kRe_k)
		\end{pmatrix},
		$$
		where $e_iRe_i$ is  a division ring, $e_iRe_i\simeq e_jRe_j$ for each $1\leq i,j\leq k$ and $n_1, \ldots, n_k $ are any positive integers.  Furthermore, if $e_iRe_j\ne 0$, then $$dim(_{e_iRe_i}(e_iRe_j))=1=dim((e_iRe_j)_{e_jRe_j}).$$
	\end{thm}

\section{Rings with every finitely generated right ideal automorphism-invariant.}

A ring $R$ is called a {\it right fa-ring (fq-ring)} if every finitely generated right ideal of $R$ is  automorphism-invariant (resp., quasi-injective).\\

\begin{rem} Let $R$ be a ring. Then, we have 
\begin{enumerate}
	\item 	  Every right  $a$-ring (or $fq$-ring) is a right  $fa$-ring.
\item If $R$ is a commutative ring,  then $R$ is a $fa$-ring if and only if it is an automorphism-invariant ring.
\item Every von Neumann regular right self-injective ring is a right $fa$-ring.
\end{enumerate}

\end{rem}

\begin{ex}\label{exam2} Consider the ring $R$ consisting of all eventually constant sequences of elements from $\mathbb F_2$ (see {\cite[Example 9]{ESS}}). Clearly, $R$ is a commutative $a$-ring. But $R$ is   not self-injective.  Thus, $R$ is a $fa$-ring and it is  not a $fq$-ring.
\end{ex}

\begin{ex}\label{exam3} The ring of linear transformations $R:= End(V_D)$ of a vector
space $V$ infinite-dimensional over a division ring $D$. It follows that $R$  is not a right $a$-ring. But $R$ is a right $fa$-ring, since every finitely generated  ideal is a direct summand of $R$ and $R$ is right self-injective.
\end{ex}

The singular submodule $Z(M)$ of a right $R$-module $M$ is defined as $Z(M)=\{m\in M: ann^{r}_R(m)$ is an essential right ideal of $R\}$ where $ann^r_R(m)$ denotes the right annihilator of $m$ in $R$.
The singular submodule of $R_R$ is called the (right) singular
ideal of the ring $R$ and is denoted by $Z(R_R)$. It is well known that $Z(R_R)$ is indeed an ideal of $R$. If $Z(R_R)=0$ then $R$ is called a  {\it right non-singular} ring.  

Recall that a ring $R$ is called {\it von Neumann regular} if for every $a\in R$, there
exists some $b\in R$ such that $a = aba$.  

\begin{prop}\cite{QAT} Let $R$ be a  regular right $fa$-ring and  $S$ be  a singular simple injective right $R$-module with  $D=\End_R(S)$. Then the  ring 
	$$
	K =
	\left(
	\begin{array}{cc}
	D & _DS_R\\
	0 & R
	\end{array}
	\right)
	$$ 
is a right $fa$-ring.
\end{prop}

\begin{cor}  A non-singular $I$-finite   ring $R$ is a right $fa$-ring if and only if it is a  semisimple artinian ring.
\end{cor}

 A ring $R$ is said to be {\it semi-primitive } if the Jacobson $J(R)=0$ and a ring $R$ is called {\it right semiartinian } if  all its nonzero factor modules of $R_R$ have nonzero socle.

We denote the ring of $n\times n$ matrices over a ring $R$ by $\mathbb M_n(R)$. In the next results  we study when matrix rings are right fa-rings.

\begin{prop}\cite{QAT}  Let  $n >1$ be an integer. The following conditions are equivalent  for a  right non-singular  ring $R$:
	\begin{enumerate}
		\item $R$ is a right self injective von Neumann regular ring.
		\item   $\mathbb M_n(R)$  is a right fa-ring.
		\item   $\mathbb M_n(R)$  is a right automorphism-invariant ring.
	\end{enumerate}
\end{prop}

\begin{cor}  Ler $R$ be a non-semisimple right self-injective von Neumann regular ring. Then for every $n >1$,  $\mathbb M_n(R)$ is a right fa-ring but it is not a right a-ring.
\end{cor}

\begin{cor}  Let $n >1$ be an integer. The following conditions are equivalent for a semi-primitive right semiartinian  ring $R$:
	\begin{enumerate}
		\item $R$ is a semisimple artinian ring.
		\item   $\mathbb M_n(R)$ is a right fa-ring.
		\item   $\mathbb M_n(R)$  is a right automorphism-invariant ring.
	\end{enumerate}
\end{cor}

In the next results  we study a decomposition  of $fa$-rings. First, recall that a right $R$-module $M$ is called {\it Zelmanowitz regular} if, every cyclic submodule of $M$ is projective and is a direct summand of $M$ (\cite{Zel}). 
  
\begin{lem}\label{key} Let $R$ be a right $fa$-ring and $A$ be a finitely generated right ideal of $R$. If there exists a right ideal $B$ of $R$ with $A\cap B=0$ and $A\cong B$, then:
\begin{enumerate}
\item $A$ is Zelmanowitz regular and injective.
\item $A$ is nonsingular.
\end{enumerate}
\end{lem}

\begin{cor} Let $e, f $ be orthogonal idempotents in a fa-ring R. If $eR\cong fR$,
then every cyclic right ideal in $eR$ is injective.
\end{cor}


Using the previous lemma, the following decomposition theorem for any right $fa$-ring is proved.

\begin{thm}\cite[Theorem 3.8]{QAT} A right $fa$-ring is isomorphic a formal matrix ring  of the form 
	$\begin{pmatrix}
	S&0\\
	M&T
	\end{pmatrix}, $	where $S$ is a square-full   von Neumann regular self-injective ring, $T$ is a  right square-free ring and $M$ is a $T-S$-bimodule.
\end{thm}

A ring $R$ is said to be {\it prime} if the product of any two nonzero ideals of $R$ is nonzero and a ring $R$ is called {\it semiprime} if it has
no nonzero nilpotent ideals.   A ring $R$ is called {\it abelian} if all idempotents of $R$ are central.

\begin{cor}\label{hqm} A semiprime  $fa$-ring is a direct sum of a square-full von Neumann regular self-injective  ring and a right square-free ring.\end{cor}

\begin{cor}\label{hqm} An abelian right $fa$-ring is a right square-free ring.
\end{cor}

\noindent A ring $R$ is said to be {\it strongly regular} if for every $a\in R$,
there exists some $b\in R$ such that $a = a^2b$.

\begin{thm}\label{at}\cite[Theorem 3.11]{QAT} The following conditions are equivalent for a right nonsingular ring $R$.

\begin{enumerate}
	\item $R$ is a right $fa$-ring.
\item $R$ is a direct sum of a square-full von Neumann regular right self-injective ring and a strongly regular ring $R$ containing
all invertible elements of its right maximal ring of fractions. 
\end{enumerate}
\end{thm}

\begin{cor}
An simple ring $R$ containing a square is a right $fa$-ring if and only if $R$ is von Neumann regular  right self-injective.
\end{cor}

\begin{thm}\cite[Theorem 3.15]{QAT}  The following conditions are equivalent for  $R$ .

\begin{enumerate}
	\item $R$ is a right nonsingular  right $fa$-ring with $\Soc(R_R)\leq^e R_R$.
    \item $R$ is a right nonsingular  right $a$-ring with $\Soc(R_R)\leq^e R_R$.
    \item $R$ is a right nonsingular  right automorphism-invariant ring with $\Soc(R_R)\leq^e R_R$.
    \item $R\cong A\times B$ where $A$ is isomorphic to a direct product of right full 
linear rings  and 
   $\oplus _{i\in I}F_i\leq R \leq \prod_{i\in I}F_i$ where $F_i\cong \mathbb{F}_2$ for every $i\in I.$  
\end{enumerate}
\end{thm}

\begin{cor}\label{cm3}  The following conditions are equivalent for right semiartinian  $R.$

\begin{enumerate}
	\item $R$ is a right nonsingular (semiprime) right $fa$-ring .
    \item $R$ is a right nonsingular (semiprime)  right $a$-ring.
    \item $R$ is a right nonsingular (semiprime)  right automorphism-invariant ring.
    \item $R\cong A\times B$ where $A$ is a  semisimple ring  and 
   $\oplus _{i\in I}F_i\leq R \leq \prod_{i\in I}F_i$ where $F_i\cong \mathbb{F}_2$ for every $i\in I.$  
\end{enumerate}
\end{cor}

The following example show that the assumption nonsingular is not superfluous in Theorem \ref{at}.

\begin{ex} 
Let $F$ be field,  $R=\prod_{i=1}^{\infty} \mathbb{M}_2(F)$, $I$ be an  essential maximal right ideal of $R$, $S=R_R/I$ and $D=\End(S)$. Consider the following ring 
	$$
	K =
	\left(
	\begin{array}{cc}
	D & S\\
	0 & R
	\end{array}
	\right).
	$$ 
	
It is obvious that	$R$ is a self-injective regular ring.  By lemma  \ref{exmn} $K$ is  a right $fa$-ring. Take  $f$ a non-trivial central idempotent of $K$. It obvious that either $f= \left(
	\begin{array}{cc}
	0 & 0\\
	0 &g
	\end{array}
	\right)$ or $f= \left(
	\begin{array}{cc}
	1_D & 0\\
	0 & 1_R-g
	\end{array}
	\right),$ where $g\in I$ is a non-trivial central idempotent of $R.$ Thus, we have a decomposition 
	$$K=\left(
	\begin{array}{cc}
	D & S\\
	0 & (1_R-g)R
	\end{array}
	\right)\oplus \left(
	\begin{array}{cc}
	D & S\\
	0 & gR
	\end{array}
	\right).$$  
	We have that  $g\neq 0$, $g\ne1$ and obtain that  rings $\left(
	\begin{array}{cc}
	D & S\\
	0 & (1_R-g)R
	\end{array}  \right)$ and  $
	\left(
	\begin{array}{cc}
	0 & 0\\
	0 & gR
	\end{array}
	\right)$ are  not square-free. Thus, the ring $K$ cannot be represented in the form $R=R_1\oplus R_2$, where $R_1$ is a regular right self-injective ring and $R_2$ is right square-free ring.
	
\end{ex}

\begin{prop}\cite{QAT}  An  indecomposable ring $R$  with nontrivial  idempotents  is a  right nonsingular right fa-ring if and only if $R$   is a von Neumann regular right self-injective ring.
\end{prop}

Recall that a ring $R$ is called \textit{directly-finite} if $xy=1$ implies $yx=1$ for all $x,y\in R$. Assume that $R$ is an abelian  right $a$-ring. By Corollary \ref{hqm}, we have a decomposition $R=S\times T$, where $S_S$ is semi-simple artinian and $T_T$ is square-free. Since $S$ and $T$ are directly-finite rings, one infers that the ring $R$ is also directly-finite.  If for a ring $R$, every matrix ring $\mathbb M_n(R)$ is directly finite then $R$ is called a {\it stably-finite ring}. It is known that the property of stable-finiteness is of importance in topology as well as in the theory of operator algebras.

\begin{prop}\cite{QAT}
Every abelian  right $fa$-ring is stably-finite.
\end{prop}

\section{Rings Whose Every Right Ideal is a Finite Direct Sum of Automorphism-Invariant Right Ideals.}

A ring $R$ is a right $\Sigma$-$a$-ring ($\Sigma$-$q$-ring, respectively) provided that its every right ideal is a finite
direct sum of automorphism-invariant right $R$-modules (quasi-injective, respectively). The $\Sigma$-$q$-rings were
introduced and studied in \cite{JSS}. Some properties of $q$-rings and $\Sigma$-$q$-rings were reflected in [22] as well.
The $\Sigma$-$a$-rings were first introduced in \cite{SS2}, and the question of description of $\Sigma$-$a$-rings was posed there
(Question 4, p. 310).

According to \cite[Theorem 2.3]{JMS}, the class of right $q$-rings coincides with the class of right self-injective
rings whose every essential right ideal is an ideal. By \cite{TQS}, a ring $R$ is a right $a$-ring if and only if $R$ is
right automorphism-invariant and each essential right ideal of R is a left $T$-module, where $T$ is a subring
of R generated by the invertible elements of $R.$ Thus, each commutative self-injective ring is a $q$-ring and
each commutative automorphism-invariant ring is an a-ring. Clearly, each $q$-ring is a $\Sigma$-$q$-ring and each
right $a$-ring is a right $\Sigma$-$a$-ring. In this section, we consider some examples of right $\Sigma$-$a$-rings showing
that the following inclusions of classes of rings are strict:

$$\{\text{ right $a$-rings} \}\subset \{\text{ right $\Sigma$-$a$-rings}\},$$
$$\{\text{ right $\Sigma$-$q$-rings} \}\subset \{\text{ right $\Sigma$-$a$-rings}\}.$$

{\bf Lemma~4.1.\ }{\sl  The following are equivalent for a quasi-Frobenius local ring $R$:
\begin{enumerate}
\item[1)]  $R$ is a right $q$-ring;
 \item[2)] $R$ is a right $a$-ring and a right $\Sigma$-$q$-ring.
 \end{enumerate}}

The following relations hold for the classes of $q$-rings, $a$-rings, $\Sigma$-$q$-rings, and  $\Sigma$-$a$-rings:  
$$
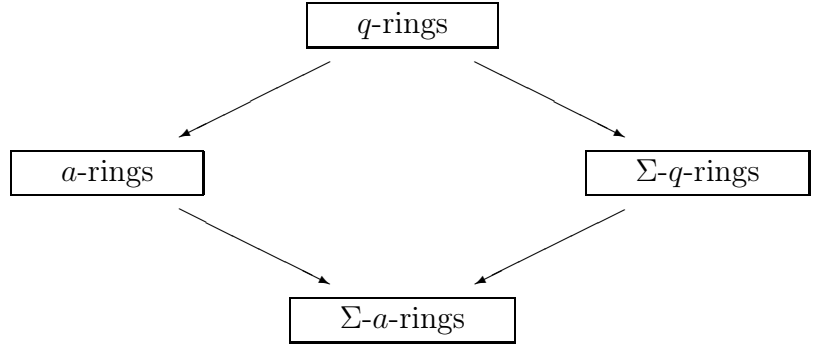
\begin{diagram}
\node[2]{\fbox{\ \ \  \ $q$-rings \ \ \ }}\arrow{sw,t}{}\arrow{se,t}{}\\
\node{\fbox{\ \ \  \  $a$-rings \ \ \ }}\arrow{se,t}{}\node[2]{\fbox{\ \ \  \  $\Sigma$-$q$-rings \ \ \ }}\arrow{sw,t}{}\\
\node[2]{\fbox{\ \ \  \  $\Sigma$-$a$-rings \ \ \ }}
\end{diagram}$$

{\bf Example.\ }  

 (a) Let $S=\prod_{i\geq 1} S_{i}$, where $\mathbb F_2=S_{i}$ for every $i$, and $$R=\{a\in S\mid \exists N:  \forall i,j>N,~a_{i}=a_{j} \}.$$
Since $E(R_R)=S_R$ and $Aut(S_R)$ is the trivial group; therefore, $R$ is an $a$-ring and, consequently, a $\Sigma$-$a$ ring. As $R$ is not self-injective, $R$ is not a $\Sigma$-q ring by Lemma 3.1.

(b) Let $F$ be a field and
$$R=\begin{pmatrix}\mathbb{F}&0&\mathbb{F}&\mathbb{F}\\ 0&\mathbb{F}&\mathbb{F}&\mathbb{F}\\ 0&0&\mathbb{F}&\mathbb{F}\\ 0&0&0&\mathbb{F}  \end{pmatrix} .$$
By \cite[Example 17]{DH}, R is a right $\Sigma$-q ring. Since $Re_{44}$ is an indecomposable and nonuniform left Rmodule, R is not left self-injective. As $R$ is an Artinian ring, $R$ is not right self-injective. Then $R$ is not
a right a-ring by Lemma 3.1.

{\bf Proposition 4.2. \cite{APT}\ }{\sl  Let $R$ be a ring, and let $S$ be the endomorphism ring  of the module $E(R_R)$. Then the
following are equivalent:
\begin{enumerate}
\item[1)] $R$ is a right $\Sigma$-$a$-ring;
\item[2)] for every right ideal $I$ of $R$ there exists a set of orthogonal idempotents $\{e_1, e_2,\dots, e_{n}\}$ in S
satisfying the conditions:
 \begin{enumerate}
        \item[a)]  $I=e_1I\oplus  e_2I\oplus\cdots\oplus e_nI$; 
        \item[b)]   $e_iI$ is a left $U_i$-module, where $U_i$  is a subring of $e_iSe_i$ generated by the invertible elements in $e_iSe_i$ for all $i=1,2,\dots,n$.
    \end{enumerate}
\end{enumerate}}

{\bf Proposition 4.3 \cite{APT}.\ }{\sl Let $e$ be an idempotent of a ring $R.$ If $R$ is a right $\Sigma$-$a$-ring and $R = ReR$ then $eRe$ is
a right $\Sigma$-$a$-ring.}

{\bf Corollary 3.4.\ }{\sl If $\mathbb{M}_n(R)$ is a right $\Sigma$-$a$-ring then $R$ is a right $\Sigma$-$a$-ring as well.}

The converse of Corollary 4.4 fails in general. To exhibit an appropriate example we need

{\bf Proposition 4.5 \cite{APT}.\ }{\sl The following are equivalent for a quasi-Frobenius local ring $R$:

(1) $R$ is a uniserial ring;

(2) $R/J(R)^2$ is a quasi-Frobenius ring;

(3) $R/J(R)^2$ is a uniserial ring;

(4) $M_2(R)$ is a right $\Sigma$-$a$-ring. }

{\bf Example.\ } 
Let $\mathbb{F}$ be a field, and let $R$ be a commutative $\mathbb{F}$-algebra with a basis
$\{1, e_0, e_1, e_2\}$  and
the product
$$1r=r1=r\ \text{for every $r\in R$,}$$
$$e_0e_1=e_1e_0=0,$$
$$e_0^2=e_1^2=e_2,$$
$$e_0e_2=e_2e_0=e_1e_2=e_2e_1=e_2^2=0.$$

 It is easy to verify that R is a quasi-Frobenius local ring and $R/J(R)^2$ is not quasi-Frobenius. Since every ideal is a quasi-injective $S$-module for an arbitrary
commutative self-injective ring $S$. Hence, every ideal of $R$ is a quasi-injective $R$-module. Then $R$ is
a $\Sigma$-$a$-ring, and $\mathbb{M}_2(R)$ is not a $\Sigma$-$a$-ring by Proposition 3.5.

A ring $R$ is right automorphism-invariant provided that $R_R$ is an automorphism-invariant module

{\bf Theorem 4.6 \cite[Theorem 7]{APT}.\ }{\sl Let $R$ be a prime right automorphism-invariant ring. The following are equivalent:

\begin{enumerate}
\item[1)] $R$ is a right $\Sigma$-$a$-ring;
\item[2)] $R$ is a Artinian ring.
\end{enumerate}}

{\bf Corollary 4.7.\ }{\sl The following conditions are equivalent for the ring of linear maps $R:=End_D(V)$ of the vector space $V$ over a division $D$.
\begin{enumerate}
\item[1)] $R$ is a right $\Sigma$-$q$-ring;
\item[2)] $R$ is a right $\Sigma$-$a$-ring;
\item[3)] the vector space $V$ is finite dimensional.
\end{enumerate}}

{\bf Corollary 4.8.\ }{\sl Let $R$ be a simple ring. Then the following conditions are equivalent:
\begin{enumerate}
\item[1)] $R$ is a right $\Sigma$-$q$-ring;
\item[2)] $R$ is a right $\Sigma$-$a$-ring;
\item[3)] $R$ is a Artinian ring.
\end{enumerate}}

{\bf Corollary 4.9.\ }{\sl Let $R$ be a primary ring. Then the following conditions are equivalent:
\begin{enumerate}
\item[1)] $R$ is a right $a$-ring;
\item[2)] $R$ is a right self-injective right  $\Sigma$-$q$-ring;
\item[3)] $R$ s a right automorphism-invariant right $\Sigma$-$a$-ring;
\item[4)] $R$ is a simple Artinian ring.
\end{enumerate}}

{\bf Proposition 4.10 \cite{APT}.\ }{\sl Let $R$ be a  right $\Sigma$-$a$-ring. If every homomorphic image of a ring $R/J(R)$ is a self-injective right-right $\Sigma$-$a$-ring, then the ring $R$ is semiperfect.}

\vskip 0.2cm 

{\bf Proposition 4.11 \cite{APT}.\ }{\sl Let $R$ be a hereditary right-Artinian right-right $\Sigma$-$a$-ring. Then the following statements hold:
\begin{enumerate}
\item[1)] if $e, f$ are two indecomposable idempotents of  $R$ and $eRf\neq 0$, then $$\dim(_{eRe}eRf)=1=\dim(eRf_{fRf});$$
\item[2)] if $R$ is an indecomposable ring, then $eRe\cong fRf$ for any nonzero indecomposable idempotents $e$ and $f$.
\end{enumerate}}

{\bf Example.\ }  
Let $X = \{1, 2, 3\}$ and $(X, \leq)$ is a partially ordered set, where $\leq=\{(1,1), (2,2), (3,3), (1,2), (1,3) \}$. Consider the incidence algebra over the field $\mathbb{F}_2$:
$$R = I(X, \mathbb{F}_2)=\begin{pmatrix}\mathbb{F}_2&\mathbb{F}_2&\mathbb{F}_2\\ 0&\mathbb{F}_2&0\\ 0&0&\mathbb{F}_2 \end{pmatrix} .$$
It is well known that $R$ is a hereditary algebra and every right $R$-module is a direct sum of the following indecomposable modules: $$e_{11}R, e_{12}R, e_{13}R, e_{11}R/e_{12}R, e_{11}R/e_{13}R, e_{11}R/(e_{12}R+e_{13}R),$$ where $e_{ij}$ are the matrix units of algebra $R$. Modules $e_{11}R/e_{12}R, e_{11}R/e_{13}R$ are injective, $e_{12}R, e_{13}R, e_{11}R/(e_{12}R+e_{13}R)$ are simple modules and $E(e_{11}R)\cong (e_{11}R/e_{12}R)\oplus (e_{11}R/e_{13}R).$ Since $\Aut(E(e_{11}R))=\{id\}$, then the module $e_{11}R$ is automorphism-invariant. Hence, every right ideal of an algebra $R$ is a direct sum of automorphism-invariant right $R$-modules. On the other hand, since the module $e_{11}R$ is not uniform, then $e_{11}R$ is not a quasi-injective module. Thus, $R$ is a hereditary right $\Sigma$-$a$-algebra that is not a right $\Sigma$-$q$-algebra.
\medskip

The following statements consider the structure of the right Artinian right hereditary right $\Sigma$-$a$-rings.

{\bf Theorem~4.12 \cite[Theorem 16]{APT}.\ }{\sl Every indecomposable right Artinian right hereditary right $\Sigma$-$a$-ring is isomorphic to
the ring of formal matrices

\[
\left[
\begin{array}{cccccc}
\mathbb{M}_{n_{1}}(D) & \mathbb{M}_{n_{1}\times n_{2}}(M_{12})
& . & . & . & \mathbb{M}_{n_{1}\times n_{m}}(M_{1m}) \\
0 & \mathbb{M}_{n_{2}}(D) & . & . & . & \mathbb{M}_{n_{2}\times
n_{m}}(M_{2m}) \\
0 & 0 & \mathbb{M}_{n_{3}}(D) & . & . & \mathbb{M}_{n_{3}\times
n_{m}}(M_{3m}) \\
. & . & . & . & . & . \\
. & . & . & . & . & . \\
0 & 0 & . & . & . & \mathbb{M}_{n_{m}}(D)
\end{array}
\right],
\]

where $D$ is a division, $M_{ij}$ be $D$-$D$ is a bimodule for every $1<i<j< m$ and $\dim(_{D}M_{ij})=1=\dim(M_{ijD})$ if $M_{ij}\neq 0$ .}

{\bf Theorem 4.13 \cite[Theorem 19]{APT}.\ }{\sl For an indecomposable right Artinian  right hereditary  ring $R$, the following conditions are equivalent:
\begin{enumerate}
\item[1)] $R$ is a serial ring; 
\item[2)] $R$ is a  right $\Sigma$-$q$-ring and  for all primitive idempotent $e$ and $f$ of  $R$ either $eRf\neq 0$ or $fRe\neq 0;$
\item[3)] $R$ is a  right $\Sigma$-$a$-ring and  for all primitive idempotent $e$ and $f$ of $R$ either $eRf\neq 0$ or $fRe\neq 0;$
\item[4)] $R$ is isomorphic to the ring of block upper triangular matrices over a skew-field.
\end{enumerate}}

{\bf Corollary~4.14 \cite[theorem 8.11]{G64}.\ }{\sl For a ring $R$, the following conditions are equivalent:
\begin{enumerate}
\item[1)] $R$ is a right non-singular  Artinian serial ring;
\item[2)]  $R$ is isomorphic to the finite direct product of the rings of block upper triangular matrices over
a skew-field.
\end{enumerate}}

\end{document}